\documentclass[10pt, twoside]{article}
\usepackage{amssymb,latexsym,amsmath}
\usepackage[hypertex]{hyperref}
\usepackage[english]{babel}
\begin{document}
\noindent
\\
\makeatletter
\renewcommand{\@evenfoot}{ \thepage \hfil \footnotesize{\it 
ISSN 1025-6415 \ \ Reports of the National Akademy of Sciences of Ukraine, 
2003, no.  7} } 
\renewcommand{\@oddfoot}{\footnotesize{\it ISSN 1025-6415 
\ \  Dopovidi Natsionalno\" \i \ Akademi\" \i \  Nauk Ukra\"\i ni, 
2003, no. 7} \hfil \thepage } 
\par{
\leftskip=1.3cm  \rightskip=0cm  
\noindent 
UDC 512 \medskip \\
\copyright \ {\bf 2003}\medskip \\
{\bf T. R. Seifullin} 
{\large \bf \medskip \\
Continuation of root functionals of a system \smallskip \\
of polynomial equations and the reduction of polynomials \smallskip \\
modulo its ideal \medskip \\ } 
{\it (Presented by Corresponding Member of the NAS of Ukraine A. A. Letichevsky)}
\medskip \\
\small {\it 
The notion of a root functional of a system of polynomials or
an ideal of polynomials is a generalization of the notion of a root 
for a multiple root. 
The operation of continuation of root functionals
and the operation of reduction of polynomials modulo the ideal are constructed
on the basis of the operation of extension of bounded root functionals
when the number of equations is equal to the number of unknowns
and the number of roots is finite.}
\par \bigskip } \noindent
The notion of a root functional arose at the investigation of linear 
relations of polynomials with polynomial coefficients  (syzygies) 
and  is a generalization of the notion of a root for the case \hbox{including} 
also and multiple roots [1-6].
An bounded root functional characterize roots of a \hbox{system} of polynomial 
equations including also and in infinity. 
A linear  functional this is an infinitely component object, 
therefore there arise the problem of its finite determination and 
operating by it in such a representation. 
An extension operation of bounded  root  functionals allows to continue 
a functional from its determination on the space  
of polynomials of the bounded degree, 
and also by this operations to reduce a polynomial modulo 
ideal, when the number of roots taking in account of multiplicity is finite.  
An extension operation is defined for a system of polynomials, 
in which the number of polynomials is equal to the number of variables. 
\smallskip

 Let ${\bf R} $ be a commutative ring with unity $1$ and zero 
$0$.

 Let $x=(x_1,\ldots ,x_n)$  be variables, ${\bf R} [x]$  be 
the ring of polynomials in variables $x$ with coefficients in ${\bf 
R} $.

 In the paper we will use definition and assumption, given 
in [6].

 Let $L(x_*)$  be a functional in ${\bf R} [x]_*$, $G(x)$  
 be a polynomial in ${\bf R} [x]$,  denote by $L(x_*)G(x)$ a functional 
with the following action:  $L(x_*)G(x).F(x)=  L(x_*).G(x)F(x)$, 
where $F(x)\in {\bf R} [x]$.

 {\bf Definition  1.}   Let  $x=(x_1,\ldots ,x_n)$,   $y\simeq 
x$    be   variables,   $f(x)   =$ $(f_1(x),\ldots ,f_n(x))$ 
 be polynomials.

 1. For a functional $L(x_*)$ denote by $\left[ L(x_*)\right] $ 
the operator \smallskip

\ \ \ $L(y_*).\det \left\| \begin{matrix}\nabla f(x,y) 
& \nabla _x(x,y) \cr
 f(x) &  {\bf 1} _x(x)\end{matrix} \right\|  = L(y_*).\det 
\left\| \begin{matrix}\nabla f(x,y) & \nabla _x(x,y) \cr
 f(y) &  {\bf 1} _x(y)\end{matrix} \right\|. $ \smallskip 

Note that since $\nabla _x(x,y)$ is linear over ${\bf R} $ operator, 
 then  operator  $\left[ L(x_*)\right] $ is linear over 
${\bf R} $.

 2.  For a functional  $L(x_*)$  and a polynomial  $F(x)$  denote by 
 $L(x_*)*F(x)=\left[ L(x_*)\right] .F(x)$,  \smallskip then

\ \ \ $L(x_*)*F(x) 
= L(y_*).\det \left\| \begin{matrix}\nabla f(x,y) & \nabla F(x,y) 
\cr
 f(x) &  F(x)\end{matrix} \right\|  = L(y_*).\det \left\| \begin{matrix}\nabla 
f(x,y) & \nabla F(x,y) \cr
 f(y) &  F(y)\end{matrix} \right\|. $ \smallskip

 3. For functionals $l(x_*)$ and $L(x_*)$ denote by $l(x_*)*L(x_*) 
= l(x_*).\left[ L(x_*)\right] $, then \smallskip

\ \ \ $l(x_*)*L(x_*) = 
  l(x_*).L(y_*).\det \left\| \begin{matrix}\nabla f(x,y) & \nabla 
_x(x,y) \cr
 f(x) &  {\bf 1} _x(x)\end{matrix} \right\| = $

\ \ \ \qquad $=\ l(x_*).  L(y_*).\det 
\left\| \begin{matrix}\nabla f(x,y) & \nabla _x(x,y) \cr
 f(y) &  {\bf 1} _x(y)\end{matrix} \right\|.$  \smallskip

This  map   $*$   we call an extension operation. 

 {\it {\bf Lemma 1.}  Let $x=(x_1,\ldots ,x_n)$  be variables, 
$f(x) =  (f_1(x),\ldots ,f_n(x))$  be polynomials, then for functionals 
$l(x_*)$, $L(x_*)$  and a polynomial  $F(x)\in {\bf R} [x]$  
there holds \smallskip

\ \ \ $l(x_*)*L(x_*).F(x) = l(x_*).L(x_*)*F(x) = l(x_*).\left[ 
L(x_*)\right] .F(x).$ } \smallskip

{\bf Proof.}\smallskip   

\ \ \ $l(x_*)*L(x_*).F(x) = 
       \left( l(x_*).\left[ L(x_*)\right] \right) .F(x) = \smallskip
       l(x_*).\left( \left[ L(x_*)\right] .F(x)\right)=$

\ \ \ \qquad $=l(x_*).L(x_*)*F(x).$ \medskip 

{\it 
 {\bf Lemma 2.}  Let $x=(x_1,\ldots ,x_n)$, $y\simeq x$  be 
variables, $f(x) =  (f_1(x),\ldots ,f_n(x))$  be polynomials, then 
for any functional  $l(x_*)$ there holds \smallskip

\ \ \ $l(x_*)*1 = l(y_*).\det \left\| \nabla f(x,y)\right\|. $ } 
\smallskip
  
{\bf Proof.}  Set $F(x)=1$, then $\nabla F(x,y) 
= 0$. We have \smallskip

\ \ \ $l(x_*)*1 = l(x_*)*F(x) =$ \smallskip

\ \ \ \qquad $ = l(y_*).\det \left\| \begin{matrix}\nabla f(x,y) & \nabla F(x,y) \cr
    f(x) &  F(x)\end{matrix} \right\| 
  = l(y_*).\det \left\| \begin{matrix}\nabla f(x,y) & 0 \cr f(x) & F(x)
    \end{matrix} \right\| =$ \smallskip

\ \ \ \qquad $ =l(y_*).\det \left\| \nabla f(x,y)\right\| \cdot F(x)      
  = l(y_*).\det \left\| \nabla f(x,y)\right\| \cdot 1   
  = l(y_*).\det \left\| \nabla f(x,y)\right\|.$ 
\medskip

{\it   
 {\bf Theorem 1.}  Let $x=(x_1,\ldots ,x_n)$  be variables, 
$f(x) =  (f_1(x),\ldots ,f_n(x))$   be polynomials, ${\delta} _f 
=\sum\limits^{ n} _{i=1} (\deg (f_i)-1)$. Let $\forall i=1,2: 
L_i(x_*)$  annuls  $(f(x))^{\leq {\delta} _f+{\delta} _i} _x$, 
where ${\delta} _i\geq 0$, then \smallskip

\ \ \ $L_1(x_*)*L_2(x_*) 
  =   L_2(x_*)*L_1(x_*)   \hbox{    in } {\bf R} [x^{\leq 
{\delta} _f+{\delta} _1+{\delta} _2+1} ].$ } \smallskip  

 {\bf Proof.}  This is the reformulation of theorem 4 in [6].  

{\it
 {\bf Theorem 2.}  Let $x=(x_1,\ldots ,x_n)$  be variables, 
$f(x) =  (f_1(x),\ldots ,f_n(x))$   be polynomials, ${\delta} _f 
=\sum\limits^{ n} _{i=1} (\deg (f_i)-1)$.

 1. Let $\forall i=1,2,3:  L_i(x_*)$  annuls  
$(f(x))^{\leq {\delta} _f+{\delta} _i} _x$,  where  ${\delta} _i\geq 0$, then 
\smallskip

\ \ \ $\left( L_1(x_*)*L_2(x_*)\right) *L_3(x_*)    =    L_1(x_*)*\left( 
L_2(x_*)*L_3(x_*)\right)     \hbox{    in } {\bf R} [x^{\leq 
{\delta} _f+{\delta} _1+{\delta} _2+{\delta} _3+2} ].$ \eject

 2.  Let  $\forall i=1,2:  L_i(x_*)$  annuls  $(f(x))^{\leq 
{\delta} _f+{\delta} _i} _x$,  where  ${\delta} _i\geq 0$,  let 
$F(x)\in {\bf R} [x^{\leq d} ]$, then \smallskip

\ \ \ $\left[ L_1(x_*)*L_2(x_*)\right] .F(x) 
\buildrel{ (f(x))^{\leq \max ({\delta} _f,d-{\delta} 
_1-{\delta} _2-2)} _x}\over \equiv  \left[ L_1(x_*)\right] .
\left[ L_2(x_*)\right] .F(x)$ \medskip
\\
and \smallskip

\ \ \ $\left[ L_2(x_*)\right] .\left[ 
L_1(x_*)\right] .F(x) \buildrel{ (f(x))^{\leq \max ({\delta} 
_f,d-{\delta} _1-{\delta} _2-2)} _x}\over \equiv  \left[ L_1(x_*)\right] 
.\left[ L_2(x_*)\right] .F(x).$ } \bigskip

{\bf Proof.}  This theorem is non-trivial and its proof 
is laborious, therefore it will be 
given in the subsequent 
papers. 
\smallskip {\it

 {\bf Theorem   3.}    Let   $x=(x_1,\ldots ,x_n)$,   $y\simeq 
x$    be   variables,   $f(x)    =$ $(f_1(x),\ldots ,f_n(x))$ 
 be  polynomials,  ${\delta} _f  =\sum\limits^{ n} _{i=1} (\deg 
(f_i)-1)$,  let  $l(x_*)$  annuls $(f(x))_x$, let $F(x)\in 
{\bf R} [x]$, then:

 1) \smallskip

\ \ \ $l(x_*)*F(x) = l(y_*)\cdot F(y).\det \left\| \nabla f(x,y)\right\| $ 
\medskip
\\
and \medskip 

\ \ \ $l(x_*)*F(x) = \left[ l(x_*)\right] .F(x) = \left( l(y_*)\cdot 
\det \| \nabla f(x,y)\| \right) .F(y);$ \bigskip

 2) $l(x_*)*F(x) \in  {\bf R} [x^{\leq {\delta} _f} ]$;

 3) $l(x_*)*F(x)$ is uniquely determined,  up to 
 addend in $(f(x))^{\leq {\delta} _f} _x$, under non-uniqueness of 
$\nabla f(x,y)$, and not depend on $\nabla F(x,y)$;

 4) if $F(x) \in  (f(x))_x$, then $l(x_*)*F(x) = 0$.   } \smallskip

 {\bf Proof 1,2.}  \smallskip \\

\ \ \ $l(x_*)*F(x) = 
      l(y_*).\det \left\| \begin{matrix}\nabla f(x,y) & \nabla F(x,y) 
      \cr f(y) &  F(y)\end{matrix} \right\|  
      = l(y_*).\det \left\| \begin{matrix}\nabla  f(x,y) & \nabla F(x,y) \cr  
      0 &  F(y)\end{matrix} \right\|=$ \bigskip

\ \ \ \qquad $=l(y_*).F(y)\cdot \det \left\| \nabla f(x,y)\right\|  
= l(y_*)\cdot F(y).\det \left\| \nabla f(x,y)\right\|.$ \bigskip
\\ 
The second equality holds, since $l(y_*)$ annuls  $(f(y))_y$. 
 The obtained polynomial $\in  {\bf R} [x^{\leq {\delta} _f} ]$, since 
$\det \left\| \nabla f(x,y)\right\| $ has a degree $\leq {\delta} 
_f$ in $x$.

 {\bf Proof  3.}  From {\it 1} of the theorem see,  that $l(x_*)*F(x)$ 
not depend  on $\nabla F(x,y)$. The functional $l(x_*)$ annuls 
$(f(x))^{\leq {\delta} _f+d} _x$ for any $d\geq 0$, then 
by virtue of {\it 2}  of theorem  2  in  [6] the polynomial 
\hbox{$l(x_*)*F(x)$} 
 is uniquely determined   up to addend belonging to 
 $(f(x))^{\leq \max ({\delta} _f,\deg (F)-d-1)} _x$,  independently 
of the choice of $\nabla f(x,y)$. For sufficiently large $d$,
 $\max ({\delta} _f,\deg (F)-d-1)  =$ ${\delta} _f$.  Hence, 
 $l(x_*)*F(x)$  is uniquely determined    up to 
addend in $(f(x))^{\leq {\delta} _f} _x$, independently of the choice of 
$\nabla f(x,y)$.

 {\bf Proof 4.}  By virtue of  {\it  1}  of the theorem there holds \medskip

\ \ \ $l(x_*)*F(x)  =  l(y_*).F(y)\cdot \det \left\| \nabla 
f(x,y)\right\|   =  0.$ \medskip
\\
The last   equality holds, 
since $F(y) \in  (f(y))_y$, and $l(y_*)$ annuls $(f(y))_y$. 
\eject {\it

 {\bf Theorem   4.}    Let   $x=(x_1,\ldots ,x_n)$,   $y\simeq 
x$     be   variables,   $f(x)   =$ $(f_1(x),\ldots ,f_n(x))$ 
 be  polynomials,  ${\delta} _f  =\sum\limits^{ n} _{i=1} (\deg 
(f_i)-1)$,  let  $l(x_*)$  annuls $(f(x))_x$, $L(x_*)$ 
annuls $(f(x))^{\leq {\delta} _f+{\delta} } _x$, where ${\delta} 
\geq 0$, then:

 1) $L(x_*)*l(x_*) = l(x_*)*L(x_*)$;

 2) $l(x_*)*L(x_*) = l(x_*)\cdot \left( L(y_*).\det \| \nabla 
f(x,y)\| \right) $;

 3) $l(x_*)*L(x_*)$ is uniquely determined, independently of the 
choice of  $\nabla f(x,y)$,  and not depend on the operator $\nabla _x(x,y)$;

 4) $l(x_*)*L(x_*)$ annuls $(f(x))_x$;

 5) $l(x_*)*L(x_*)$ not depend on the action of $L(x_*)$ outside ${\bf 
R} [x^{\leq {\delta} _f} ]$. }

 {\bf Proof 1.}  Since $l(x_*)$ annuls $(f(x))^{\leq 
{\delta} _f+d} _x$ for any  $d\geq 0$, $L(x_*)$ annuls 
$(f(x))^{\leq {\delta} _f+{\delta} } _x$, then by virtue of theorem 1 
$l(x_*)*L(x_*) = L(x_*)*l(x_*)$ in 
${\bf R} [x^{\leq {\delta} _f+{\delta} +d+1} ]$, 
and mean, and in the whole ${\bf R} [x]$  
by the arbitrariness of $d\geq 0$.

 {\bf Proof 2.}  Let $F(x)\in {\bf R} [x]$, then \medskip

\ \ \ $l(x_*)*L(x_*).F(x) = 
   l(x_*).L(y_*).\det \left\| \begin{matrix}\nabla f(x,y) & \nabla F(x,y) \cr
 f(x) &  F(x)\end{matrix} \right\|  =$ \medskip
 
\ \ \ $\qquad = l(x_*).L(y_*).\det \left\| \begin{matrix}\nabla f(x,y) & 
  \nabla F(x,y) \cr 0 &  F(x)\end{matrix} \right\|
  = l(x_*).L(y_*).\det \left\| \nabla f(x,y)\right\| \cdot F(x) = $ \medskip
 
\ \ \ $\qquad = l(x_*)\cdot \left( L(y_*).
      \det \| \nabla f(x,y)\| \right) .F(x).$ \bigskip
\\
From the arbitrariness of $F(x)$ we have the equality of functionals \medskip

\ \ \ $l(x_*)*L(x_*) = l(x_*)\cdot \left( L(y_*).\det \| \nabla f(x,y)\| \right). $ 
\medskip

 {\bf Proof 3.}  From  {\it  2}  see,   that  $l(x_*)*L(x_*)$ 
 not  depend  on $\nabla _x(x,y)$.

 Since $l(x_*)$ annuls $(f(x))^{\leq {\delta} _f+d} 
_x$ for any $d\geq 0$,  $L(x_*)$  annuls $(f(x))^{\leq 
{\delta} _f+{\delta} } _x$, then by virtue of {\it 1}  of theorem  3  in 
 [6]  $l(x_*)*L(x_*)$  is uniquely determined  in ${\bf R} 
[x^{\leq {\delta} _f+{\delta} +d+1} ]$  independently  of the   choice of 
  $\nabla f(x,y)$,   and   mean, is uniquely determined in the whole 
${\bf R} [x]$  by the arbitrariness of $d\geq 0$.

 {\bf Proof 4.}  From {\it 2} see,  that  $l(x_*)*L(x_*)$ 
 annuls  $(f(x))_x$, since $l(x_*)$ annuls $(f(x))_x$.

 {\bf Proof 5.}  From  {\it  2}  see,   that  $l(x_*)*L(x_*)$ 
 not  depend  on the action of $L(x_*)$ outside ${\bf R} [x^{\leq {\delta} 
_f} ]$, since $\det \left\| \nabla f(x,y)\right\| $ has a degree  
$\leq {\delta} _f$  in  $y$. \smallskip

 {\bf Definition 2.} Let $x=(x_1,\ldots ,x_n)$, $y\simeq 
x$ \   be   variables, let $f(x)= (f_1(x),\ldots ,f_n(x))$  be polynomials.  
A functional  $E(x_*)$  we call a unit  root functional  of polynomials  
$f(x)$,  if  it  annuls  $(f(x))_x$,  and  
$E(x_*)*1   =E(y_*).\det \left\| \nabla f(x,y)\right\| 
 =  1  +  f(x)\cdot g(x)$. A functional  $E'(x_*)$  we call a  unit 
bounded   root  functional  of polynomials  $f(x)$,  if 
  it   annuls $(f(x))^{\leq {\delta} _f+{\varepsilon} } _x$, 
where  ${\varepsilon} \geq 0$,  and  $E'(x_*)*1 = E'(y_*).\det 
\left\| \nabla f(x,y)\right\|  = 1 + f(x)\cdot g(x)$. \smallskip  {\it 

 {\bf Lemma 3.}  Let $x=(x_1,\ldots ,x_n)$, $y\simeq x$  be 
variables, $f(x) =  (f_1(x),\ldots ,f_n(x))$  be polynomials. Let 
a functional  $E(x_*)$  annuls  $(f(x))_x$,  and  $E(x_*)*1  
-  1  \in $ $(f(x))_x$, let $E'(x_*) = E(x_*)$ in ${\bf R} 
[x^{\leq {\delta} _f+{\varepsilon} } ]$, where ${\varepsilon} \geq 
0$. Then  $E'(x_*)$  annuls 
$(f(x))^{\leq {\delta} _f+{\varepsilon} } _x$ and 
$E'(x_*)*1 -1 \in  (f(x))_x$. } \eject

 {\bf Proof.}  $E'(x_*)*1 =  E'(y_*).\det \left\| 
\nabla f(x,y)\right\|   =  E(y_*).\det \left\| \nabla f(x,y)\right\| 
  =$ $E(x_*)*1$, since $E'(y_*) = E(y_*)$ in ${\bf R} [y^{\leq 
{\delta} _f+{\varepsilon} } ]$ and  $\det \left\| \nabla f(x,y)\right\| 
$  has a degree $\leq {\delta} _f$ in $y$. Hence, 
$E'(x_*)*1 -1 =$ \hbox{$E(x_*)*1$} $-1 \in  (f(x))_x$. Since  $E'(x_*)$ $= E(x_*)$ 
in ${\bf R} [x^{\leq {\delta} _f+{\varepsilon} } ]$ and $E(x_*)$ 
annuls $(f(x))^{\leq {\delta} _f+{\varepsilon} } _x \subseteq 
 {\bf R} [x^{\leq {\delta} _f+{\varepsilon} } ]$, then  and  $E'(x_*)$ 
annuls $(f(x))^{\leq {\delta} _f+{\varepsilon} } _x$. \smallskip 
 {\it

 {\bf Theorem   5.}    Let   $x=(x_1,\ldots ,x_n)$,   $y\simeq 
x$    be   variables,   $f(x)    =$ $(f_1(x),\ldots ,f_n(x))$ 
 be polynomials,  ${\delta} _f  =\sum\limits^{ n} _{i=1} (\deg (f_i)-1)$. 
 Let a  functional  $E'(x_*)$ annuls $(f(x))^{\leq {\delta} 
_f+{\varepsilon} } _x$, where ${\varepsilon} \geq 0$, and $E'(x_*)*1 
- 1 \in  (f(x))_x$, then: \medskip 

 1) if $F(x) \in  {\bf R} [x^{\leq d} ]$, then \bigskip

\ \ \ $E'(x_*)*F(x) \in  {\bf R} [x^{\leq \max ({\delta} _f,d-{\varepsilon} -1)} ]$
\smallskip \\
and \smallskip

\ \ \ $F(x) - E'(x_*)*F(x) \in  (f(x))_x\cap {\bf R} [x^{\leq 
\max ({\delta} _f,d)} ];$\bigskip 

 2) if $l(x_*)$ annuls $(f(x))_x$, then \bigskip

\ \ \ $l(x_*) = E'(x_*)*l(x_*) = l(x_*)*E'(x_*);$\bigskip 

 3) if $l(x_*)$ annuls $(f(x))_x$ and $L'(x_*) = l(x_*)$ 
in ${\bf R} [x^{\leq {\delta} _f+{\delta} } ]$, where ${\delta} 
\geq 0$, then \bigskip

\ \ \ $l(x_*) = E'(x_*)*L'(x_*) = L'(x_*)*E'(x_*) \hbox{ 
in } {\bf R} [x^{\leq {\delta} _f+{\delta} +{\varepsilon} +1} 
];$ \bigskip

 4) if a functional $L(x_*)$ annuls $(f(x))_x\cap {\bf 
R} [x^{\leq {\delta} _f+{\delta} } ]$, where  ${\delta} \geq 0$, 
then \bigskip

\ \ \ $L(x_*)  =  E'(x_*)*L(x_*)  =  L(x_*)*E'(x_*) \hbox{ in } 
{\bf R} [x^{\leq {\delta} _f+{\delta} } ]$\bigskip
\\
and $E'(x_*)*L(x_*)=L(x_*)*E'(x_*)$ 
annuls $(f(x))_x\cap {\bf R} [x^{\leq {\delta} _f+{\delta} 
+{\varepsilon} +1} ]$. }

 {\bf Proof 1.}   Since  $E'(x_*)$  annuls 
 $(f(x))^{\leq {\delta} _f+{\varepsilon} } _x$  and  
$F(x)  \in $ ${\bf R} [x^{\leq d} ]$, then by virtue of {\it 1} of 
theorem 2 in  [6] 
 $E'(x_*)*F(x)  \in   {\bf R} [x^{\leq \max ({\delta} _f,d-{\varepsilon} 
-1)} ]$. Then \bigskip

\ \ \ $F(x) - E'(x_*)*F(x) 
\in  {\bf R} [x^{\leq d} ] + 
{\bf R} [x^{\leq \max ({\delta} _f,d-{\varepsilon} -1)} ] =$\bigskip

\ \ \ \qquad $= {\bf R} [x^{\leq \max (d,\max ({\delta} _f,d-{\varepsilon} -1))} ] 
 =   {\bf R} [x^{\leq \max ({\delta} _f,d)} ]. $ \medskip
\\ 
Moreover,  \medskip

\ \ \ $E'(x_*)*F(x) =$ \medskip

\ \ \ $\qquad = E'(y_*).\det \left\| \begin{matrix}\nabla f(x,y) 
& \nabla F(x,y) \cr
 f(x) &  F(x)\end{matrix} \right\|    \buildrel{ (f(x))_x}\over 
\equiv E'(y_*).\det \left\| \begin{matrix}\nabla f(x,y) & \nabla 
F(x,y) \cr  0 &  F(x)\end{matrix} \right\| =$ \medskip

\ \ \ $\qquad =  \left( E'(y_*).\det \| \nabla f(x,y)\| \right) 
\cdot F(x)  =  \left( 1+f(x)\cdot g(x)\right) \cdot F(x)   \buildrel{ 
(f(x))_x}\over \equiv F(x),$ \bigskip
\\
then   $F(x)   
-   E'(x_*)*F(x)    \in     (f(x))_x$.    
Finally,
$F(x)-E'(x_*)*F(x)\in 
(f(x))_x\cap {\bf R} [x^{\leq \max ({\delta} _f,d)} ]$.

 {\bf Proof 2.}  Since $l(x_*)$ annuls $(f(x))_x$, 
 $E'(x_*)$  annuls $(f(x))^{\leq {\delta} _f+{\varepsilon} 
} _x$,  then  by virtue of  {\it  1}  of theorem  4 
there holds $E'(x_*)*l(x_*)  = l(x_*)*E'(x_*)$, 
and by virtue of {\it 2} of theorem 4 \medskip

\ \ \ $l(x_*)*E'(x_*) = 
l(x_*)\cdot (E'(y_*).\det \left\| \nabla f(x,y)\right\| ) = 
l(x_*)\cdot (1+f(x)\cdot g(x)) = l(x_*).$ \medskip
\eject\noindent
The last equality 
 holds, since  $l(x_*)$ annuls $(f(x))_x$, and since $l(x_*)\cdot 1 = l(x_*)$.

 {\bf Proof 3.}  Since $L'(x_*)$ \  annuls 
 $(f(x))^{\leq {\delta} _f+{\delta} } _x$,  and  $E'(x_*)$ annuls 
 $(f(x))^{\leq {\delta} _f+{\varepsilon} } _x$, 
 then by virtue of theorem 1 $L'(x_*)*E'(x_*) = E'(x_*)*L'(x_*)$ 
 in ${\bf R} [x^{\leq {\delta} _f+{\delta} +{\varepsilon} +1} 
]$. Since and $L'(x_*) = l(x_*)$ in ${\bf R} [x^{\leq {\delta} 
_f+{\delta} } ]$,  then by virtue of {\it 3} of theorem 3 in [6] $l(x_*)*E'(x_*) 
= L'(x_*)*E'(x_*)$ in ${\bf R} [x^{\leq {\delta} _f+{\delta} 
+{\varepsilon} +1} ]$. By virtue of {\it 2}  of the theorem  $l(x_*)  =  l(x_*)*E'(x_*)$. 
 Hence, there holds $l(x_*)=L'(x_*)*E'(x_*) = E'(x_*)*L'(x_*)$ 
in ${\bf R} [x^{\leq {\delta} _f+{\delta} +{\varepsilon} +1} 
]$.

 {\bf Proof 4.}  Let $F(x) \in   {\bf R} [x^{\leq 
{\delta} _f+{\delta} } ]$,  then  $\max ({\delta} _f,{\delta} 
_f+{\delta} )  =  {\delta} _f+{\delta} $, since ${\delta} 
\geq 0$. By virtue of the second statement of {\it 1} of the theorem \medskip

\ \ \ $F(x)- E'(x_*)*F(x) \in  (f(x))_x\cap {\bf R} [x^{\leq \max ({\delta} 
_f,{\delta} _f+{\delta} )} ] = (f(x))_x\cap {\bf R} [x^{\leq 
{\delta} _f+{\delta} } ],$
\medskip
\\
and therefore is annulled by $L(x_*)$. 
We have, by using of the second equality of lemma 1, \medskip

\ \ \ $0 = L(x_*).(F(x) - E'(x_*)*F(x)) = L(x_*).F(x) - L(x_*)*E'(x_*).F(x),$
\medskip
\\ 
hence,  from the arbitrariness of \ 
$F(x) \in  {\bf R} [x^{\leq {\delta} _f+{\delta} } ]$, 
$L(x_*) = L(x_*)*E'(x_*)$ \ in ${\bf R} [x^{\leq {\delta} _f+{\delta} } ]$.  
Since  $E'(x_*)$   annuls   
$(f(x))^{\leq {\delta} _f+{\varepsilon} } _x$,   
$L(x_*)$   annuls 
$(f(x))_x\cap {\bf R} [x^{\leq {\delta} _f+{\delta} } ] \supseteq  
(f(x))^{\leq {\delta} _f+{\delta} }_x$, 
and mean, and 
$(f(x))^{\leq {\delta} _f+{\delta} } _x$, then by virtue of theorem  1 
$E'(x_*)*L(x_*) = L(x_*)*E'(x_*)$ in 
${\bf R} [x^{\leq {\delta} _f+{\delta} +{\varepsilon} +1} ] \supseteq  
{\bf R} [x^{\leq {\delta} _f+{\delta} } ]$, 
and  hence, and 
in ${\bf R} [x^{\leq {\delta} _f+{\delta} } ]$.

 Let  $F(x)  \in   (f(x))^{\leq d} _x\cap {\bf R} [x^{\leq 
{\delta} _f+{\delta} +{\varepsilon} +1} ]$.   Since   $E'(x_*)$ 
  annuls $(f(x))^{\leq {\delta} _f+{\varepsilon} } _x$,  
 then   by virtue of the first  statement of  {\it 1}   $E'(x_*)*F(x) 
   \in $ ${\bf R} [x^{\leq \max ({\delta} _f,{\delta} _f+({\delta} 
+{\varepsilon} +1)-{\varepsilon} -1)} ] = {\bf R} [x^{\leq \max 
({\delta} _f,{\delta} _f+{\delta} )} ] =  {\bf R} [x^{\leq {\delta} 
_f+{\delta} } ]$, and by virtue of {\it 3} of theorem  2  in  [6]  $E'(x_*)*F(x) 
 \in   (f(x))^{\leq d-{\varepsilon} -1} _x$,  hence, 
the polynomial $E'(x_*)*F(x) \in   (f(x))_x\cap {\bf R} [x^{\leq {\delta} 
_f+{\delta} } ]$.  From  the last  it follows that there holds 
the equality  $L(x_*)*E'(x_*).F(x)  =  L(x_*).E'(x_*)*F(x)  =   
0$,   since   $L(x_*)$ annuls $(f(x))_x\cap {\bf R} [x^{\leq 
{\delta} _f+{\delta} } ]$. Hence,  from the arbitrariness of 
a polynomial $F(x)   \in     (f(x))_x\cap {\bf R} [x^{\leq {\delta} 
_f+{\delta} +{\varepsilon} +1} ]$, the functional  $L(x_*)*E'(x_*)$ 
   annuls $(f(x))_x\cap {\bf R} [x^{\leq {\delta} _f+{\delta} 
+{\varepsilon} +1} ]$. \smallskip  {\it

 {\bf Theorem 6.}    Let   $x=(x_1,\ldots ,x_n)$,   $y\simeq 
x$    be   variables,   $f(x)    =$ $(f_1(x),\ldots ,f_n(x))$ 
 be  polynomials,  ${\delta} _f  =\sum\limits^{ n} _{i=1} (\deg 
(f_i)-1)$.  Let  $E(x_*)$  annuls $(f(x))_x$ and $E(x_*)*1 
- 1 \in  (f(x))_x$, then:

 1) if $F(x) \in  {\bf R} [x^{\leq d} ]$, then \medskip

\ \ \ $E(x_*)*F(x) 
\in  {\bf R} [x^{\leq {\delta} _f} ]$ \medskip
\\
and \medskip

\ \ \ $F(x) - E(x_*)*F(x) 
\in  (f(x))_x\cap {\bf R} [x^{\leq \max ({\delta} _f,d)} ];$
\medskip

 2) if $l(x_*)$ annuls $(f(x))_x$, then \medskip

\ \ \ $l(x_*) = E(x_*)*l(x_*) = l(x_*)*E(x_*) = 
E(x_*)\cdot \left( l(y_*).\det \| \nabla f(x,y)\| \right) ;$
\medskip

 3) if $l(x_*)$ annuls $(f(x))_x$, and $L'(x_*) = l(x_*)$ 
in ${\bf R} [x^{\leq {\delta} _f+{\delta} } ]$, where ${\delta} 
\geq 0$, then \medskip

\ \ \ $l(x_*) = E(x_*)*L'(x_*) = L'(x_*)*E(x_*);$ \medskip

 4) if  a functional  $L(x_*)$  annuls  
$(f(x))_x\cap {\bf R} [x^{\leq {\delta} _f+{\delta} } ]$,  
where  ${\delta} \geq 0$,  then \medskip

\ \ \ $L(x_*) = E(x_*)*L(x_*) = L(x_*)*E(x_*) 
\hbox{ in }  {\bf R} [x^{\leq {\delta} _f+{\delta} } ];$
\medskip
  
note  that $E(x_*)*L(x_*) = L(x_*)*E(x_*)$ annuls $(f(x))_x$; \eject

 5) if $l(x_*)$ annuls $(f(x))_x$, then  it  is uniquely determined its the  
action on ${\bf R} [x^{\leq {\delta} _f+{\delta} } ]$.  }

 {\bf Proof.}  For any ${\varepsilon} \geq 0$ the functional 
$E(x_*)$ annuls $(f(x))^{\leq {\delta} _f+{\varepsilon} } 
_x$.

 {\bf Proof 1.}  From {\it 1} of theorem 5 it follows that \medskip

\ \ \ $F(x) - E(x_*)*F(x) \in (f(x))_x\cap {\bf R} 
[x^{\leq \max ({\delta} _f,d)} ]$
\medskip
\\
and $\forall {\varepsilon} \geq 0:$ $E(x_*)*F(x)\in {\bf R} 
[x^{\leq \max ({\delta} _f,d-{\varepsilon} -1)} ]$, 
hence, $E(x_*)*F(x) \in  {\bf R} [x^{\leq {\delta} _f} ]$.

 {\bf Proof 2.}  From {\it 2} of theorem 5 it follows that \medskip

\ \ \ $l(x_*) = E(x_*)*l(x_*) = l(x_*)*E(x_*),$
\medskip
\\
and since $E(x_*)$ and $l(x_*)$ annul $(f(x))_x$,  
then from {\it 2} of theorem 4 it follows that \medskip

\ \ \ $l(x_*)*E(x_*)  = E(x_*)\cdot \left( l(y_*).\det 
\| \nabla f(x,y)\| \right) .$
\medskip

 {\bf Proof 3 and 5.} From {\it 3} of theorem 5 it follows that \medskip

\ \ \ $l(x_*) = E(x_*)*L'(x_*) = L'(x_*)*E(x_*) \hbox{ in }  
{\bf R} [x^{\leq {\delta} _f+{\delta} +{\varepsilon} +1} ].$  
\medskip
\\
From the arbitrariness of ${\varepsilon} \geq 0$ we obtain that \medskip

\ \ \ $l(x_*) = E(x_*)*L'(x_*) = L'(x_*)*E(x_*) \hbox{ in } {\bf R} [x].$
\medskip
\\
Since the equality  holds for  any  $L'(x)$  such 
 that $L'(x_*) = l(x_*)$ in ${\bf R} [x^{\leq {\delta} _f+{\delta} 
} ]$, then $l(x_*)$ is uniquely determined its  the action in ${\bf R} [x^{\leq {\delta} _f} ]$.

 {\bf Proof 4.}   The first  statement it follows  from 
 the first  statement of {\it 4} of theorem 5. The second statement it follows 
from {\it 1} and {\it 4} of theorem 4, since $L(x_*)$ annuls 
$(f(x))^{\leq {\delta} _f+{\delta} } _x$ and $E(x_*)$ annuls 
$(f(x))_x$. \medskip {\it

 {\bf Theorem   7.}    Let   $x=(x_1,\ldots ,x_n)$,   $y\simeq 
x$     be   variables,   $f(x)   =$ $(f_1(x),\ldots ,f_n(x))$ 
 be polynomials, \  ${\delta} _f  =\sum\limits^{ n} _{i=1} (\deg (f_i)-1)$. 
\  Let \  $E'(x_*)$  annuls $(f(x))_x$, \ or annuls $(f(x))^{\leq 
{\delta} _f+{\varepsilon} } _x$, where ${\varepsilon} \geq 0$, 
 and  $E'(x_*)*1  -  1  \in   (f(x))_x$. Then:

 1)   ${\bf R} [x]/(f(x))_x$   coincide   with   the set of  
 all   elements    of the form $E'(x_*)*F(x)/(f(x))_x$, where $F(x) 
\in  {\bf R} [x^{\leq {\delta} _f} ]$. Moreover,  $E'(x_*)*F(x)\in 
{\bf R} [x^{\leq {\delta} _f} ]$.

 2) $(f(x))_x\cap {\bf R} [x^{\leq d'} ]$, \ where $d'\geq {\delta} _f$, 
\ coincide  with  the set  of all  elements of the form 
\hbox{$F(x)-E'(x_*)*F(x)$}, 
where $F(x) \in  {\bf R} [x^{\leq d'} ]$.  }

 {\bf Proof.}  If the functional $E'(x_*)$ annuls 
$(f(x))_x$, then for  any ${\varepsilon} \geq 0$ it annuls 
$(f(x))^{\leq {\delta} _f+{\varepsilon} } _x$. Therefore the statement 
 it suffices to prove  for the case, when the functional $E'(x_*)$ 
annuls $(f(x))^{\leq {\delta} _f+{\varepsilon} } _x$ for 
${\varepsilon} \geq 0$.

 Consider the sequence of polynomials   $G_0(x)$,   $\forall 
p\geq 0:   G_{p+1} (x)   =$ $E'(x_*)*G_p(x)$.

 {\bf Proof 1.}  Let  $G_0(x)  \in   {\bf R} [x^{\leq 
d} ]$,  since  $E'(x_*)$  annuls $(f(x))^{\leq {\delta} 
_f+{\varepsilon} } _x$, by virtue of the second statement of 
{\it 1} of theorem 5 there hold \medskip

\ \ \ $G_0(x)-G_1(x) = G_0(x) - E'(x_*)*G_0(x) 
\in  (f(x))_x,$ \medskip

\ \ \ $\ldots \ldots \ldots \ldots \ldots \ldots \ldots \ldots \ldots 
\ldots \ldots \ldots \ldots \ldots \ldots \ldots \ldots \ldots\ldots ,$\medskip

\ \ \ $G_{p-1} (x)-G_p(x) = G_{p-1} (x) - E'(x_*)*G_{p-1} (x) \in 
 (f(x))_x.$ \eject\noindent
Hence, $S(x) = G_0(x) - E'(x_*)*G_{p-1} (x) = G_0(x) - G_p(x) \in (f(x))_x$,     
then 
there holds 
$G_0(x)=$ \hbox{$E'(x_*)*G_{p-1} (x) + S(x)$}. 
And by virtue of the first statement of {\it 1} of theorem 
 5 there hold \medskip

\ \ \ $G_1(x) = E'(x_*)*G_0(x) \in  
{\bf R} [x^{\leq \max ({\delta} _f,d-({\varepsilon} +1))} ],$

\ \ \ $\ldots \ldots \ldots \ldots \ldots \ldots \ldots \ldots \ldots 
\ldots \ldots \ldots \ldots \ldots \ldots\ldots . \,. ,$ \medskip

\ \ \ $G_p(x)  =  E'(x_*)*G_{p-1} (x)  \in    {\bf R} [x^{\leq 
\max ({\delta} _f,d-p\cdot ({\varepsilon} +1))} ].$ \medskip
\\
For sufficiently large $p$ there holds 
$\max ({\delta} _f,d-(p-1)\cdot ({\varepsilon} +1)) = {\delta} _f$.  
Hence, 
$F(x) = G_{p-1} (x) \in  {\bf R} [x^{\leq {\delta} _f} ]$. 
We finally obtain,  that any polynomial $G_0(x)$ is of the form $G_0(x) = E'(x_*)*F(x) 
+ S(x)$, where $S(x) \in  (f(x))_x$, $F(x) \in  {\bf R} [x^{\leq 
{\delta} _f} ]$. Hence, any polynomial in ${\bf R} [x]/(f(x))_x$ 
 is of the form  $E'(x_*)*F(x)/(f(x))_x$,  where  $F(x)  \in $ ${\bf 
R} [x^{\leq {\delta} _f} ]$.

 Otherwise, if $F(x)\in {\bf R} [x^{\leq {\delta} _f} ]$, 
 then $E'(x_*)*F(x)/(f(x))_x \in {\bf R} [x]/(f(x))_x$. 
Moreover, since $E'(x_*)$  annuls  $(f(x))^{\leq {\delta} 
_f+{\varepsilon} } _x$, then  by virtue of the first statement of {\it 1} 
of theorem 5 there holds $E'(x_*)*F(x)  \in   {\bf R} [x^{\leq 
\max ({\delta} _f,{\delta} _f-{\varepsilon} -1)} ]  =$ ${\bf 
R} [x^{\leq {\delta} _f} ]$.

 {\bf Proof 2.}  Let $G_0(x) \in  (f(x))_x\cap {\bf 
R} [x^{\leq d'} ]$, then $G_0(x) \in   (f(x))^{\leq d} _x$ 
for some $d$. 
Since $E'(x_*)$ annuls $(f(x))^{\leq {\delta} _f+{\varepsilon} 
} _x$, by virtue of {\it 3} of theorem 2 in [6] there hold \medskip

\ \ \ $G_1(x) = E'(x_*)*G_0(x) \in  (f(x))^{\leq d-({\varepsilon} +1)} _x,$
\medskip

\ \ \ $\ldots \ldots \ldots \ldots \ldots \ldots \ldots \ldots \ldots 
\ldots \ldots \ldots \ldots \ldots.\,.\, ,$ \medskip 

\ \ \ $G_p(x) = E'(x_*)*G_{p-1} (x) \in 
(f(x))^{\leq d-p\cdot ({\varepsilon} +1)} _x. $\medskip 

For sufficiently large $p: d-p\cdot ({\varepsilon} 
+1) < 0$, hence, $G_p(x)  =  0$.  Set  $P=p$, then 
$$G_0(x) = \sum\limits^{P-1} _{p=0} \left( G_p(x)-E'(x_*)*G_p(x)\right) 
 = \left( \sum\limits^{P-1} _{p=0} G_p(x)\right) -E'(x_*)*\left( 
\sum\limits^{P-1} _{p=0} G_p(x)\right). \hphantom{cc}$$ 
By virtue of the first statement of {\it 1} of theorem 5 
$\forall p\geq 0: G_p(x) \in  
{\bf R} [x^{\leq \max ({\delta} _f,d'-p\cdot ({\varepsilon} +1))} ]  \subseteq 
$ ${\bf R} [x^{\leq d'} ]$, since $d'\geq {\delta} _f$. 
 Hence,  any  
$G_0(x)  \in $ $(f(x))_x\cap 
{\bf R} [x^{\leq d'} ]$ \ is of the form $F(x)-E'(x_*)*F(x)$, 
 where  $F(x)  =  \sum\limits^{P-1} _{p=0} G_p(x)  \in $ ${\bf 
R} [x^{\leq d'} ]$.

 Otherwise, let $F(x)\in {\bf R} [x^{\leq d'} ]$; 
then by virtue of the second  statement of  {\it  1} of theorem 5 
the polynomial $F(x) -  E'(x_*)*F(x)$  
belongs to  $(f(x))_x\cap {\bf R} [x^{\leq \max ({\delta} _f,d')} ]= 
(f(x))_x\cap {\bf R} [x^{\leq d'} ]$, since $d'\geq {\delta} _f$. \bigskip 

{\footnotesize

\begin{enumerate}

\item {\it Seifullin, T. R.} 
Root functionals and root polynomials 
of a system of polynomials. (Russian)
Dopov. Nats. Akad. Nauk Ukra\"\i ni  -- 1995, -- no. 5, 5--8.
\item {\it Seifullin, T. R.} Root functionals and root relations 
of a system of polynomials. (Russian) 
Dopov. Nats. Akad. Nauk Ukra\"\i ni  -- 1995, -- no 6, 7--10.
\item  {\it Seifullin, T. R.}  Homology of the Koszul complex of a 
system of polynomial equations. (Russian)
Dopov. Nats. Akad. Nauk Ukr. Mat. Prirodozn. Tekh. Nauki 1997, no. 9, 43--49. 
\item  {\it Seifullin, T. R.}  Koszul complexes of systems of 
polynomials connected by linear dependence. (Russian) 
Some problems in 
contemporary mathematics (Russian), 326--349, Pr. Inst. Mat. Nats. Akad. Nauk 
Ukr. Mat. Zastos., 25, Natsional. Akad. Nauk Ukra\"\i ni, Inst. Mat., Kiev, 
1998.
\item {\it Seifullin, T. R.} Koszul complexes of embedded systems of 
polynomials and duality. (Russian) 
Dopov. Nats. Akad. Nauk Ukr. Mat. Prirodozn. 
Tekh. Nauki 2000, no. 6, 26--34. \eject  
\item {\it Seifullin, T. R.} 
Extension of bounded root functionals 
of a system of polynomial equations. Dopov. Nats. Akad. Nauk Ukr. Mat. 
Prirodozn. Tekh. Nauki 2002, no. 7, 35--42. 
\href{http://arxiv.org/abs/0804.2420}
{{\tt  arxiv:0804.2420}}.
\item {\it  Buchberger B.} Gr\"obner: An algorithmic method in polynomial
ideal theory
//Multidimensional Systems Theory. / Ed. N. K. Bose, --
Dordrecht: D. Reidel, 1985. -- Chapter 6.
\item {\it Caniglia L., Galligo A., Heintz J.} Some new effictivity bounds
in computational geometry.
// Proc. 6th Int. Conf. on Appied Algebra and Error--correcting codes.
/ LNCS 357, Springer--Verlag, Berlin. -- 1989. -- pp. 131--152.
\item {\it Brownawell D.} Bounds for the degrees in the Nullstellensatz. 
// Ann. Math. 2nd series. -- 1987. -- No 126. -- pp. 577--591.
\item {\it Canny J.} Generalized characteristic polynomials. 
//\ J.Symbolic Computation. -- 1990. -- No 9. -- pp. 241--250.
\\
\end{enumerate}

\small{\noindent
{\it V. M. Glushkov Institute of Cybernetics of the NAS of Ukraine, Kiev
\hfill Received 26.06.2002 \medskip\\
E-mail: \ {\tt  timur\_sf@mail.ru}
} 

\end{document}